\newcommand{\xBc}{\langle}
\newcommand{\xBe}{\rangle}

\newcommand{\xbD}{\Delta}
\newcommand{\xbF}{\Phi}
\newcommand{\xbG}{\Gamma}

\newcommand{\xbP}{\Pi}

\newcommand{\xbS}{\Sigma}

\newcommand{\xba}{\alpha}
\newcommand{\xbb}{\beta}

\newcommand{\xbe}{\in}
\newcommand{\xbf}{\phi}
\newcommand{\xbg}{\gamma}

\newcommand{\xbm}{\mu}

\newcommand{\xbq}{\psi}
\newcommand{\xbr}{\rho}
\newcommand{\xbs}{\sigma}
\newcommand{\xbt}{\tau}

\newcommand{\xCI}{{\Big(}}
\newcommand{\xCJ}{{\Big)}}
\newcommand{\xCK}{\times}
\newcommand{\xCL}{\pm}

\newcommand{\xCN}{\neg}
\newcommand{\xCO}{ }
\newcommand{\xCQ}{\emptyset}

\newcommand{\xCf}{\hspace{0.1em}}

\newcommand{\xcA}{\forall}

\newcommand{\xcE}{\exists}

\newcommand{\xcM}{\not\models}
\newcommand{\xcN}{\hspace{0.2em}\not\sim\hspace{-0.9em}\mid\hspace{0.8em}}

\newcommand{\xcU}{\bigwedge}

\newcommand{\xcc}{\subseteq}
\newcommand{\xcd}{\supseteq}
\newcommand{\xce}{\not\in}

\newcommand{\xcg}{\geq}
\newcommand{\xch}{\Rightarrow}
\newcommand{\xci}{\Leftarrow}
\newcommand{\xcj}{\Leftrightarrow}
\newcommand{\xck}{\leq}
\newcommand{\xcl}{\vdash}
\newcommand{\xcm}{\models}
\newcommand{\xcn}{\hspace{0.2em}\sim\hspace{-0.9em}\mid\hspace{0.58em}}

\newcommand{\xco}{\vee}
\newcommand{\xcp}{\rightarrow}

\newcommand{\xcs}{\cap}
\newcommand{\xcu}{\wedge}
\newcommand{\xcv}{\cup}

\newcommand{\xcz}{\Box}

\newcommand{\xDH}{\item }

\newcommand{\xDO}{\circ}

\newcommand{\xdf}{{\cal F}}

\newcommand{\xdj}{{\cal J}}

\newcommand{\xEI}{\begin{itemize}}
\newcommand{\xEJ}{\end{itemize}}

\newcommand{\xET}{\%}

\newcommand{\xEd}{\neq}

\newcommand{\xEh}{\begin{enumerate}}

\newcommand{\xEj}{\end{enumerate}}

\newcommand{\xeb}{\prec}
\newcommand{\xec}{\preceq}

\newcommand{\xex}{\upharpoonright}

\newcommand{\xfA}{\mid}

\newcommand{\Xl}{\ldots}

\newcommand{\bl}{\begin{lemma} \rm}
\newcommand{\el}{\end{lemma}}
\newcommand{\br}{\begin{remark} \rm}
\newcommand{\er}{\end{remark}}
\newcommand{\be}{\begin{example} \rm}
\newcommand{\ee}{\end{example}}
\newcommand{\bco}{\begin{corollary} \rm}
\newcommand{\eco}{\end{corollary}}
\newcommand{\bc}{\begin{claim} \rm}
\newcommand{\ec}{\end{claim}}
\newcommand{\bfa}{\begin{fact} \rm}
\newcommand{\efa}{\end{fact}}
\newcommand{\bp}{\begin{proposition} \rm}
\newcommand{\ep}{\end{proposition}}
\newcommand{\bd}{\begin{definition} \rm}
\newcommand{\ed}{\end{definition}}
\newcommand{\bcs}{\begin{construction} \rm}
\newcommand{\ecs}{\end{construction}}
\newcommand{\bcd}{\begin{condition} \rm}
\newcommand{\ecd}{\end{condition}}
\newcommand{\bt}{\begin{theorem} \rm}
\newcommand{\et}{\end{theorem}}
\newcommand{\bn}{\begin{notation} \rm}
\newcommand{\en}{\end{notation}}
\newcommand{\bfi}{\begin{bild} \rm}
\newcommand{\efi}{\end{bild}}
\newcommand{\bsta}{\begin{statement} \rm}
\newcommand{\esta}{\end{statement}}
\newcommand{\bcom}{\begin{comment} \rm}
\newcommand{\ecom}{\end{comment}}
\newcommand{\bdia}{\begin{diagram} \rm}
\newcommand{\edia}{\end{diagram}}

\newcommand{\bfc}{\begin{figure}[htb] \begin{center}}
\newcommand{\efc}{\end{center} \end{figure}}

\documentclass{article}

\usepackage{amssymb,latexsym,epic,eepic,rotating}

\title{
Semantic interpolation
\thanks{
paper 337
}
}

\author{Dov M Gabbay
\thanks{
Dov.Gabbay@kcl.ac.uk, www.dcs.kcl.ac.uk/staff/dg
} \\
King's College, London
\thanks{
Department of Computer Science, King's College London, Strand,
London WC2R 2LS, UK
} \\
and \\
Bar-Ilan University, Israel
\thanks{
Department of Computer Science,
Bar-Ilan University,
52900 Ramat-Gan, Israel
} \\
and \\
University of Luxembourg
\thanks{
Computer Science and Communications,
Faculty of Sciences,
6, rue Coudenhove-Kalergi,
L-1359 Luxembourg
} \\ \\
Karl Schlechta
\thanks{
ks@cmi.univ-mrs.fr, karl.schlechta@web.de, http://www.cmi.univ-mrs.fr/ $\sim$ ks
} \\
Laboratoire d'Informatique Fondamentale de Marseille
\thanks{
UMR 6166, CNRS and Universit\'{e} de Provence,
Address: CMI, 39, rue Joliot-Curie, F-13453 Marseille Cedex 13, France
}
}

\begin{document}

\newtheorem{lemma}{Lemma}[section]
\newtheorem{theorem}[lemma]{Theorem}
\newtheorem{proposition}[lemma]{Proposition}
\newtheorem{corollary}[lemma]{Corollary}
\newtheorem{claim}[lemma]{Claim}
\newtheorem{fact}[lemma]{Fact}
\newtheorem{remark}[lemma]{Remark}
\newtheorem{definition}{Definition}[section]
\newtheorem{construction}{Construction}[section]
\newtheorem{condition}{Condition}[section]
\newtheorem{example}{Example}[section]
\newtheorem{notation}{Notation}[section]
\newtheorem{bild}{Figure}[section]
\newtheorem{comment}{Comment}[section]
\newtheorem{statement}{Statement}[section]
\newtheorem{diagram}{Diagram}[section]

\renewcommand{\labelenumi}
  {(\arabic{enumi})}
\renewcommand{\labelenumii}
  {(\arabic{enumi}.\arabic{enumii})}
\renewcommand{\labelenumiii}
  {(\arabic{enumi}.\arabic{enumii}.\arabic{enumiii})}
\renewcommand{\labelenumiv}
  {(\arabic{enumi}.\arabic{enumii}.\arabic{enumiii}.\arabic{enumiv})}

\maketitle

\setcounter{secnumdepth}{3}
\setcounter{tocdepth}{3}

\begin{abstract}

We treat interpolation for various logics. In the full non-monotonic case, we
connect the existence
of interpolants to laws about abstract size, and to generalized Hamming
relations. We also mention revision a la Parikh.

\end{abstract}

\tableofcontents

%
%
%
\section{
Introduction
}

The interpolation property is the ability to ``squeeze'' a new formula
in the consequence relation where the new formula has only symbols
common to both original formulas:

If $ \xba \xcn_{1} \xbg,$ then there is an interpolant $ \xbb $ s.t. $
\xba \xcn_{2} \xbb \xcn_{3} \xbg $ - and $ \xbb $
has only symbols common to $ \xba $ and $ \xbg.$

Here, the consequence relations $ \xcn_{i}$ may be different. In the
classical
case, they are all the same, and equal to $ \xcl.$ We will see later that
the
more general case is important in other logics.

Interpolation for classical logic is well-known,
see  \cite{Cra57}, and there are also
non-classical logics for which interpolation has been shown, e.g., for
Circumscription, see
 \cite{Ami02}. An extensive overview of interpolation is found
in  \cite{GM05}. Chapter 1 of this book gives a survey and a
discussion and
11 points of view of interpolation are discussed. This paper is a
new 12th point of view.

Omitting a variable in classical propositional logic means
``liberating'' it, i.e., semantically, any of the two values
TRUE and FALSE may be taken. Thus, for the omitted variables, we
take the full product. Thus, there is a straightforward translation
from language to sequences of values TRUE, FALSE (for each variable),
and ``liberated'' variables correspond to taking all possible
combinations. We work immediately with such sequences. This has the
advantage that we need not fix our approach to just 2 values, TRUE and
FALSE, but can work with many truth values. We are also free from
logically equivalent re-formulations, classical consequence becomes
subset inclusion, etc., in short: life becomes simpler, and new
perspectives open up naturally.

We will see that changing logic, e.g. from preferential logic to
classical logic when ``squeezing'' $ \xbb $ in will be important.

Most of these notes are purely semantical, or, better, even algebraic.
In the background is (partly classical) propositional logic, but we tried
to remain general in our approach.

Our approach has the advantage of short and elementary proofs.
Perhaps more important, it separates algebraic from logical
questions, and we see that there are logics with algebraic
interpolation, but without logical interpolation, as the
necessary sets of models are not definable in the language.
This opens the way to making the language richer to
obtain interpolation, when so desired.
\section{
Notation
}

\bd

$\hspace{0.01em}$


\label{Definition Sin-Restriction}

Let $ \xbP:= \xbP \{X_{i}:i \xbe J\}.$

For $J' \xcc J$ and $ \xbs \xbe \xbP,$ let $ \xbs \xex J' $ be the
restriction of $ \xbs $ to $J',$
likewise for $ \xbS \xcc \xbP $ and $ \xbS \xex J'.$

\ed

\bd

$\hspace{0.01em}$


\label{Definition Sin-Ir-Relevant}

For $ \xbS \xcc \xbP $ set

$I(\xbS):=\{i \xbe J:$ $ \xbS = \xbS \xex (J-\{i\}) \xCK X_{i}\}$ (up to
re-ordering), (the irrelevant
or inessential $i)$ and

$R(\xbS):=J-I(\xbS)$ (the relevant or essential $i).$

\ed

\bfa

$\hspace{0.01em}$


\label{Fact Sin-Ir-Relevant}

(1) $ \xbS = \xbS \xex R(\xbS) \xCK \xbP \xex I(\xbS)$ (up to
re-ordering)

(2) $ \xbs \xex R(\xbS)= \xbs' \xex R(\xbS)$ $ \xcu $ $ \xbs \xbe
\xbS $ $ \xch $ $ \xbs' \xbe \xbS.$

\efa

\subparagraph{
Proof
}

$\hspace{0.01em}$


(1) Enumerate $I(\xbS)$ somehow. Take the first $i \xbe I(\xbS)$ where
it fails. But $i$
is irrelevant, contradiction.

$ \xcz $
\\[3ex]

\bn

$\hspace{0.01em}$


\label{Notation Products}

We consider products over $X,$ let $X',X'' $ be a disjoint cover of $X.$
Let
(by abuse of language)
$ \xbP:= \xbP X,$ $ \xbP':= \xbP X',$ etc. $ \xbs $ will be some
element of $ \xbP,$ $ \xbs' $ the restriction
of $ \xbs $ to $X',$ etc. $ \xbS $ will be a set of such $ \xbs,$ etc.
\section{
The semantically monotonic case (upward or downward)
}

\en

Note that the following interpolation results can be read upward
(monotonic logic) or downward (non-monotonic logic, in the following
sense:
$ \xbg $ is the theory of the minimal models of $ \xba,$ and not just any
formula
which holds in the set of minimal models - which would be downward,
and then upward again in the sense of set inclusion), in the latter case
we have to be careful: we usually cannot go upward again, so we have the
sharpest possible case in mind. The case of mixed movement - down and then
up -
as in full non-monotonic logic is treated in Section \ref{Section Non-Mon-Int}
(page \pageref{Section Non-Mon-Int}).

Moreover, the logic can be many-valued. There is no restriction on the
$X_{i}.$

\bfa

$\hspace{0.01em}$


\label{Fact Sin-Interpolation}

Let $ \xbS' \xcc \xbS \xcc \xbP,$ where $ \xbP = \xbP \{X_{i}:i \xbe
J\}.$

(1) Let $ \xbS'':= \xbS' \xex (R(\xbS) \xcs R(\xbS')) \xCK \xbP
\xex (I(\xbS) \xcv I(\xbS')).$

Then $ \xbS' \xcc \xbS'' \xcc \xbS.$

The following two results
concern ``parallel interpolation'', terminology introduced by
D.Makinson in  \cite{KM07}. Thus in the first case, $ \xbS' $ is a
product, in the second
case, $ \xbS $ is a product.

(2) Let $ \xdj $ be a disjoint cover of $J.$

Let $ \xbS' = \xbP \{ \xbS'_{K}:K \xbe \xdj \}$ with $ \xbS'_{K} \xcc
\xbP \{X_{i}:i \xbe K\}.$

Let $ \xbS''_{K}:= \xbS'_{K} \xex (R(\xbS) \xcs R(\xbS'_{K})) \xCK
\xbP \{X_{i}:i \xbe K,i \xbe I(\xbS) \xcv I(\xbS'_{K})\}.$

Let $ \xbS'':= \xbP \{ \xbS''_{K}:K \xbe \xdj \}$ (re-ordered).

Then $ \xbS' \xcc \xbS'' \xcc \xbS.$

(3) Let $ \xdj $ be a disjoint cover of $J.$

Let $ \xbS = \xbP \{ \xbS_{K}:K \xbe \xdj \}$ with $ \xbS_{K} \xcc \xbP
\{X_{i}:i \xbe K\}.$

Let $ \xbS''_{K}:= \xbS' \xex (R(\xbS') \xcs R(\xbS_{K})) \xCK \xbP
\{X_{i}:i \xbe K,i \xbe I(\xbS_{K}) \xcv I(\xbS')\},$
so $ \xbS''_{K} \xcc \xbP \{X_{i}:i \xbe K\}.$

Let $ \xbS'':= \xbP \{ \xbS''_{K}:K \xbe \xdj \}$ (re-ordered).

Then $ \xbS' \xcc \xbS'' \xcc \xbS.$

\efa

\subparagraph{
Proof
}

$\hspace{0.01em}$


(1)

(1.1) $ \xbS' \xcc \xbS'' $ is trivial.

(1.2) $ \xbS'' \xcc \xbS:$

Let $ \xbs \xbe \xbS'',$ so there is $ \xbs' \xbe \xbS' $ s.t.
$ \xbs \xex (R(\xbS) \xcs R(\xbS'))= \xbs' \xex (R(\xbS) \xcs R(
\xbS')).$ By definition of $R(\xbS')$
(or Fact \ref{Fact Sin-Ir-Relevant} (page \pageref{Fact Sin-Ir-Relevant})  (2)
for $ \xbS'),$
we may choose $ \xbs' \xbe \xbS' $ s.t. also $ \xbs \xex I(\xbS')=
\xbs' \xex I(\xbS'),$ so there is
$ \xbs' \xbe \xbS' \xcc \xbS $ s.t. $ \xbs \xex R(\xbS)= \xbs' \xex
R(\xbS),$ so by definition of $R(\xbS)$
(or Fact \ref{Fact Sin-Ir-Relevant} (page \pageref{Fact Sin-Ir-Relevant})  (2)
for $ \xbS),$
$ \xbs \xbe \xbS.$

(2)

(2.1) $ \xbS' \xcc \xbS''.$

$ \xbS'_{K} \xcc \xbS''_{K},$ so by $ \xbS'' = \xbP \xbS''_{K}$ the
result follows.

(2.2) $ \xbS'' \xcc \xbS.$

Let $ \xbs'' \xbe \xbS'',$ $ \xbs'' = \xDO \{ \xbs''_{K}:K \xbe \xdj
\}$ (concatenation) for suitable $ \xbs''_{K} \xbe \xbS''_{K}.$
Consider $ \xbs''_{K}.$ By definition of $ \xbS''_{K},$ there is $ \xbs
'_{K} \xbe \xbS'_{K}$ s.t.
$ \xbs''_{K} \xex (R(\xbS'_{K}) \xcs R(\xbS))= \xbs'_{K} \xex (R(
\xbS'_{K}) \xcs R(\xbS)),$ so there is $ \xbt'_{K} \xbe \xbS'_{K}$
s.t.
$ \xbs''_{K} \xex R(\xbS)= \xbt'_{K} \xex R(\xbS).$
Let $ \xbt':= \xDO \{ \xbt'_{K}:K \xbe \xdj \},$ so by $ \xbS' = \xbP
\{ \xbS'_{K}:K \xbe \xdj \}$ $ \xbt' \xbe \xbS' \xcc \xbS.$
But $ \xbs'' \xex R(\xbS)= \xbt' \xex R(\xbS),$ so $ \xbs'' \xbe
\xbS.$

(3)

We first show $R(\xbS_{K})=R(\xbS) \xcs K.$

Let $i \xbe I(\xbS_{K}),$ then $ \xbS_{K}= \xbS_{K} \xex (K-\{i\}) \xCK
X_{i},$ but $ \xbS = \xbP \{ \xbS_{K}:K \xbe \xdj \},$ so
$ \xbS = \xbS \xex (J-\{i\}) \xCK X_{i},$ and $i \xbe I(\xbS).$
Conversely, let $i \xbe I(\xbS) \xcs K,$ then
$ \xbS = \xbS \xex (J-\{i\}) \xCK X_{i},$ so $ \xbS \xex K= \xbS \xex
(K-\{i\}) \xCK X_{i},$ so $i \xbe I(\xbS_{K}).$

(3.1) $ \xbS' \xcc \xbS''.$

$ \xbS' \xex K \xcc \xbS''_{K},$ so by $ \xbS'' = \xbP \{ \xbS''_{K}:K
\xbe \xdj \}$ $ \xbS' \xcc \xbS''.$

(3.2) $ \xbS'' \xcc \xbS.$

By $ \xbS = \xbP \{ \xbS_{K}:K \xbe \xdj \},$ it suffices to show $ \xbS
''_{K} \xcc \xbS_{K}.$

Let $ \xbs''_{K} \xbe \xbS''_{K}.$ So there is $ \xbs' \xbe \xbS' $
s.t.
$ \xbs' \xex (R(\xbS_{K}) \xcs R(\xbS'))= \xbs''_{K} \xex (R(
\xbS_{K}) \xcs R(\xbS')),$ so there is $ \xbt' \xbe \xbS' \xcc \xbS $
s.t.
$ \xbt' \xex R(\xbS_{K})= \xbs''_{K} \xex R(\xbS_{K}),$ so there is $
\xbs \xbe \xbS $ s.t. $ \xbs \xex K= \xbs''_{K},$ so $ \xbs''_{K} \xbe
\xbS_{K}.$

$ \xcz $
\\[3ex]

We have shown $ \xCf semantic$ interpolation, this is not yet
$ \xCf syntactic$ interpolation. We still need that the set of
sequences is definable.
(The importance of definability in the context of non-monotonic
logics was examined by one of the authors
in  \cite{Sch92}.) Note that we
``simplified'' the set of sequences, but perhaps the logic
at hand does not share this perspective. Consider, e.g.,
intuitionistic logic with three worlds, linearly ordered.
This is a monotonic logic, so by our result, it has
semantic interpolation. But it has no syntactic
interpolation, so the created set of models must not be
definable. In classical propositional logic, the created set
$ \xCf is$ definable, as we will see in
Fact \ref{Fact Sin-Simplification-Definable} (page \pageref{Fact
Sin-Simplification-Definable}).

We describe now above example in detail.

\be

$\hspace{0.01em}$


\label{Example Sin-4-Value}

The semantics is a set of 3 worlds, $w,$ $w',$ $w'',$ linearly ordered.
The logic is intuitionistic, so knowledge
``grows''. Thus, $ \xbf $ can be true from $w$ on, from $w' $ on, from
$w'' $ on, or never, resulting in 4 truth values.

The formulas to consider are

$ \xba (p,q,r):= \xCI p \xcp (((q \xcp r) \xcp q) \xcp q) \xCJ \xcp p,$

$ \xbb (p,s):=((s \xcp p) \xcp s) \xcp s$

We have provable $ \xba \xcp \xbb,$ but no syntactic interpolant (which
could use only $p).$

Introducing a new operator Jp meaning
``from next moment onwards $p$ holds and if now is the last moment then $p$
holds
now''
gives enough definability to have also syntactic interpolation.

\ee

We may also see this definability property as natural, and as a
criterion for well-behaviour of a logic. In some cases, introducing
new constants analogous to TRUE, FALSE - in the cited case e.g. ONE, TWO
when truth starts at world one or two - might help, but we did not
investigate this question.
This question is also examined in  \cite{ABM03}.

We can go further with a logic in language L0 for which there is no
interpolation. For every pair of formulas which give a counterexample to
interpolation we introduce a new connective which corresponds to the
semantic
interpolant. Now we have a language L1 which allows for interpolation
for formulas in the original language L0. L1 itself may or may not
have interpolation. So we might have to continue to L2, L3, etc.

This leads to the following definition for the cases where we have
semantic, but not necessarily syntactic interpolation:

\bd

$\hspace{0.01em}$


\label{Definition Sin-Level}

A logic has level 0 semantic interpolation, iff it has interpolation.

A logic has level $n+1$ semantic interpolation iff it has not level $n$
semantic interpolation, but introducing new elements into the language
(of level $n)$ results in interpolation also for the new
language.

\ed

We have not examined this notion.

The case of full non-monotonic logic is, of course, different, as
the logics might not even have semantic interpolation,
so above repairing is not possible.

\bfa

$\hspace{0.01em}$


\label{Fact Sin-Simplification-Definable}

Simplification preserves definability in classical propositional logic:

Let $ \xbG = \xbS \xex X' \xCK \xbP X''.$
Then, if $ \xbS $ is formula definable, so is $ \xbG.$

\efa

\subparagraph{
Proof
}

$\hspace{0.01em}$


As $ \xbS $ is formula definable, it is defined by $ \xbf_{1} \xco  \Xl
\xco \xbf_{n},$ where
$ \xbf_{i}= \xbq_{i,1} \xcu  \Xl  \xcu \xbq_{i,n_{i}}.$ Let $ \xbF_{i}:=\{
\xbq_{i,1}, \Xl, \xbq_{i,n_{i}}\},$
$ \xbF'_{i}:=\{ \xbq \xbe \xbF_{i}: \xbq \xbe X' \}$ (more precisely, $
\xbq $ or $ \xCN \xbq \xbe X'),$ $ \xbF''_{i}:= \xbF_{i}- \xbF'_{i}.$
Let $ \xbf'_{i}:= \xcU \xbF'_{i}.$ Thus $ \xbf_{i} \xcl \xbf'_{i}.$ We
show that $ \xbf'_{1} \xco  \Xl  \xco \xbf'_{n}$ defines $ \xbG.$
(Alternatively, we may replace all $ \xbq \xbe \xbF''_{i}$ by TRUE.)

(1) $ \xbG \xcm \xbf'_{1} \xco  \Xl  \xco \xbf'_{n}$

Lt $ \xbs \xbe \xbG,$ then there is $ \xbt \xbe \xbS $ s.t. $ \xbs \xex
X' = \xbf \xex X'.$ By prerequisite,
$ \xbt \xcm \xbf_{1} \xco  \Xl  \xco \xbf_{n},$ so $ \xbt \xcm \xbf'_{1}
\xco  \Xl  \xco \xbf'_{n},$ so $ \xbs \xcm \xbf'_{1} \xco  \Xl  \xco
\xbf'_{n}.$

(2) Suppose $ \xbs \xce \xbG,$ we have to show $ \xbs \xcM \xbf'_{1}
\xco  \Xl  \xco \xbf'_{n}.$

Suppose then $ \xbs \xce \xbG,$ but $ \xbs \xcm \xbf'_{1} \xco  \Xl
\xco \xbf'_{n},$ wlog. $ \xbs \xcm \xbf'_{1}= \xcU \xbF'_{1}.$
As $ \xbs \xce \xbG,$ there is no $ \xbt \xbe \xbS $ $ \xbt \xex X' =
\xbs \xex X'.$
Choose $ \xbt $ s.t. $ \xbs \xex X' = \xbt \xex X' $ and $ \xbt \xcm \xbq
$ for all $ \xbq \xbe \xbF''_{i}.$
By $ \xbs \xcm \xbq $ for $ \xbq \xbe \xbF'_{1},$ and $ \xbs \xex X' =
\xbt \xex X' $ $ \xbt \xcm \xbq $ for $ \xbq \xbe \xbF'_{1}.$
By prerequisite, $ \xbt \xcm \xbq $ for $ \xbq \xbe \xbF''_{1},$ so $
\xbt \xcm \xbq $ for all $ \xbq \xbe \xbF_{1},$ so
$ \xbt \xcm \xbf_{1} \xco  \Xl  \xco \xbf_{n},$ and $ \xbt \xbe \xbS,$ as
$ \xbf \xco  \Xl  \xco \xbf $ defines $ \xbS,$ contradiction.

$ \xcz $
\\[3ex]

\bco

$\hspace{0.01em}$


\label{Corollary Sin-Simplification-Definable}

The same result holds if $ \xbS $ is theory definable.

\eco

\subparagraph{
Proof
}

$\hspace{0.01em}$


(Outline). Define $ \xbS $ by a (possibly infinite) conjunction of
(finite)
disjunctions. Transform this into a possibly infinite disjunction of
possibly infinite conjunctions. Replace all $ \xbf \xbe \xbF''_{i}$ by
TRUE. The same
proof as above shows that this defines $ \xbG $ (finiteness was nowhere
needed). Transform backward into a conjunction of finite disjunctions,
where
the $ \xbf \xbe \xbF''_{i}$ are replaced by TRUE.

$ \xcz $
\\[3ex]
\section{
The full non-monotonic case, i.e., downward and upward
}

\label{Section Non-Mon-Int}
\subsection{
Discussion
}

We consider here a non-monotonic logic $ \xcn.$
$ \xcn $ will be defined by a principal filter, generated by a set
operator $ \xbm (X) \xcc X:$

$ \xbb \xcn \xbg $ iff $ \xbm (M(\xbb)) \xcc M(\xbg),$ where $M(\xba
)$ is the set of models of $ \xba.$

Often, $ \xbm $ will be generated by a binary relation $ \xeb $ on the
model set,

$ \xbm (X):=\{x \xbe X: \xCN \xcE x' \xeb x,$ $x' \xbe X\}.$

The following example shows that interpolation does not always exist for
the full non-monotonic case.

\be

$\hspace{0.01em}$


\label{Example Sin-Non-Mon-Int}

Full non-monotonic logics, i.e. down and up, has not necessarily
interpolation.

Consider the model order $pq \xeb p \xCN q \xeb \xCN p \xCN q \xeb \xCN
pq.$ Then $ \xCN p \xcn \xCN q,$ there are no
common variables, and $true \xcn q$ (and, of course, $ \xCN p \xcN
false).$
(Full consequence of $ \xCN p$ is $ \xCN p \xCN q,$ so this has
trivial interpolation.)

\ee

We look at the interpolation problem in 3 ways.

Given $ \xbf \xcn \xbq,$ there is an interpolant $ \xba $ s.t.

(1) $ \xbf \xcn \xba \xcl \xbq $

(2) $ \xbf \xcl \xba \xcn \xbq $

(3) $ \xbf \xcn \xba \xcn \xbq $

The first variant will be fully characterized below.

The second variant has no full characterization at the time of writing
(to the authors' knowledge), but is connected to very interesting
properties about multiplication of size and componentwise independent
relations.

The third variant stays unexplored for the moment.

We turn to variant (1) and (2).
\subsection{
$\xbf \xcn \xba \xcl \xbq $
}

\label{Section Non-Mon-Int-Karl}

\bfa

$\hspace{0.01em}$


\label{Fact Sin-Int-Aux}

Let $ \xbS \xcc \xbP X,$ $var(\xba) \xcs var(\xbb)= \xCQ,$ $var(\xbb
) \xcs R(\xbS)= \xCQ,$ $ \xbb $ not a tautology, then
$ \xbS \xcc M(\xba \xco \xbb)$ $ \xch $ $ \xbS \xcc M(\xba).$

\efa

\subparagraph{
Proof
}

$\hspace{0.01em}$


Suppose not, so there is $ \xbs \xbe \xbS $ such that $ \xbs \xcm \xba
\xco \xbb,$ $ \xbs \xcM \xba.$ As $ \xbb $ is not
a tautology, there is an assignment to $var(\xbb)$ which makes $ \xbb $
wrong. Consider
$ \xbt $ s.t. $ \xbs = \xbt $ except on $var(\xbb),$ where $ \xbt $
makes $ \xbb $ wrong, using this assignment.
By $var(\xba) \xcs var(\xbb)= \xCQ,$ $ \xbt \xcm \xCN \xba.$ By
$var(\xbb) \xcs R(\xbS)= \xCQ,$ $ \xbt \xbe \xbS.$
So $ \xbt \xcM \xba \xco \xbb $ for some $ \xbt \xbe \xbS,$
contradiction.

$ \xcz $
\\[3ex]

\bfa

$\hspace{0.01em}$


\label{Fact Sin-Non-Mon-Int-Karl}

We use here normal forms (conjunctions of disjunctions).

Consider a finite language.

$ \xcn $ has interpolation iff for all $ \xbS,$ $ \xbm (\xbS)$ $I(\xbS
) \xcc I(\xbm (\xbS))$ holds.

In the infinite case, we need as additional prerequisite that
$ \xbm (\xbS)$ is definable if $ \xbS $ is.

$(\xbm (\xbS)$ are the minimal models of $ \xbS,$ seen as sequences.)

\efa

\subparagraph{
Proof
}

$\hspace{0.01em}$


Work with reformulations of $ \xbS $ etc. which use only essential $(=$
relevant)
variables.

``$ \xch $'':

Suppose the condition is wrong. Then $X:=I(\xbS)-I(\xbm (\xbS))=I(
\xbS) \xcs R(\xbm (\xbS)) \xEd \xCQ.$
Then there is some $ \xbs' \xbe \xbm (\xbS) \xex R(\xbS)$ which
cannot be continued by some choice $ \xbr $ in $X \xcv (I(\xbS) \xcs I(
\xbm (\xbS)))$ in $ \xbm (\xbS),$ i.e.
$ \xbs' \xDO \xbr \xce \xbm (\xbS).$

We first do the finite case: Consider the
formula $ \xbf:= \xbs' \xcp \xCN \xbr = \xCN \xbs' \xco \xCN \xbr.$ We
have $Th(\xbS) \xcn \xbf,$ as $ \xbm (\xbS) \xcc M(\xbf).$
Suppose $ \xbS'' $ is a semantical interpolant for $ \xbS $ and $ \xbf.$
So $ \xbm (\xbS) \xcc \xbS'' \xcc M(\xbf),$
and $ \xbS'' $ does not contain any variables in $ \xbr $ as essential
variables.
By Fact \ref{Fact Sin-Int-Aux} (page \pageref{Fact Sin-Int-Aux}), $ \xbm (\xbS
) \xcc \xbS'' \xcc M(
\xCN \xbs'),$ but $ \xbs' \xbe \xbm (\xbS) \xex R(\xbS),$
contradiction.

We turn to the infinite case. Consider again $ \xbs' \xDO \xbr.$ As $
\xbs' \xDO \xbr \xce \xbm (\xbS),$
and $ \xbm (\xbS)$ is definable,
there is some formula $ \xbf $ which holds in $ \xbm (\xbS),$ but fails
in $ \xbs' \xDO \xbr.$
Thus, $Th(\xbS) \xcn \xbf.$ Write $ \xbf $ as a disjunction of
conjunctions.
Let $ \xbS'' $ be an interpolant for $ \xbS $ and $M(\xbf).$ Thus $
\xbm (\xbS) \xcc \xbS'' \xcc M(\xbf),$ and
$ \xbs' \xDO \xbr \xce M(\xbf),$ so $ \xbm (\xbS) \xcc \xbS'' \xcc
M(\xbf) \xcc M(\xCN \xbs' \xco \xCN \xbr),$ so
$ \xbS'' \xcc M(\xCN \xbs')$ by Fact \ref{Fact Sin-Int-Aux} (page
\pageref{Fact Sin-Int-Aux}).
So $ \xbm (\xbS) \xcm \xCN \xbs',$ contradiction, as $ \xbs' \xbe
\xbm (\xbS) \xex R(\xbS).$

``$ \xci $'':

Let $I(\xbS) \xcc I(\xbm (\xbS)).$ Let $ \xbS \xcn \xbS',$ i.e. $
\xbm (\xbS) \xcc \xbS'.$ Write $ \xbm (\xbS)$ as a (possibly
infinite) conjunction of
disjunctions, using only relevant variables. Form $ \xbS'' $ from $ \xbm
(\xbS)$
by omitting all variables in this description which are not in $R(\xbS'
).$
Note that all remaining variables are in $R(\xbm (\xbS)) \xcc R(\xbS
),$ so
$ \xbS'' $ is a candidate for interpolation.

(1) $ \xbm (\xbS) \xcc \xbS'':$ Trivial.

(2) $ \xbS'' \xcc \xbS':$
Let $ \xbs \xbe \xbS''.$ Then there is $ \xbt \xbe \xbm (\xbS) \xcc
\xbS' $ s.t. $ \xbs \xex R(\xbS')= \xbt \xex R(\xbS'),$ so
$ \xbs \xbe \xbS'.$

It remains to show in the infinite case that $ \xbS'' $ is definable.
This can be
shown as in Fact \ref{Fact Sin-Simplification-Definable} (page \pageref{Fact
Sin-Simplification-Definable}).

$ \xcz $
\\[3ex]
\subsection{
$\xbf \xcl \xba \xcn \xbq $
}
\subsubsection{
Product size
}

We work here with a notion of ``big'' and
``small'' subsets, which may be thought of as defined by a filter (ideal),
though we usually will not need the full strength of a filter (ideal). But
assume as usual that $A \xcc B \xcc C$ and $A \xcc C$ is big imply $B \xcc
C$ is big, that $C \xcc C$ is
big, and define
$A \xcc B$ is small iff $(B-A) \xcc B$ is big, call all subsets which are
neither big
nor small medium size. For an extensive discussion, see
 \cite{GS09a}.

Let $X' \xcv X'' =X$ be a disjoint cover, so $ \xbP X= \xbP X' \xCK \xbP
X''.$ We consider subsets $ \xbS $
etc. of $ \xbP X.$ If not said otherwise, $ \xbS $ etc. need not be a
product $ \xbS' \xCK \xbS''.$
As usual, $ \xbP':= \xbP X',$ $ \xbS'':= \xbS \xex X'' $ etc. The
roles of $X' $ and $X'' $ are
interchangeable, e.g., instead of $ \xbG \xex X' \xcc \xbS \xex X',$ we
may also write
$ \xbG \xex X'' \xcc \xbS \xex X''.$

We consider here the following three finite product rules about size
and $ \xbm.$ $ \xbm (X) \xcc X$ - read: the minimal elements of $X$ -
will generate
the principal filter $ \xdf (X),$ whose elements are big subsets of $X.$

\bd

$\hspace{0.01em}$


\label{Definition Sin-Size-Rules}

$(S*1)$ $ \xbD \xcc \xbS' \xCK \xbS'' $ is big iff there is $ \xbG =
\xbG' \xCK \xbG'' \xcc \xbD $ s.t. $ \xbG' \xcc \xbS' $ and
$ \xbG'' \xcc \xbS'' $ are big

$(S*2)$ $ \xbG \xcc \xbS $ is big $ \xch $ $ \xbG \xex X' \xcc \xbS \xex
X' $ is big - where
$ \xbS $ is not necessarily a product.

$(S*3)$ $A \xcc \xbS $ is big $ \xch $ there is $B \xcc \xbP' \xCK \xbS
'' $ big s.t. $B \xex X'' \xcc A \xex X'' $ - again,
$ \xbS $ is not necessarily a product.

$(\xbm *1)$ $ \xbm (\xbS' \xCK \xbS'')= \xbm (\xbS') \xCK \xbm (
\xbS'')$

$(\xbm *2)$ $ \xbm (\xbS) \xcc \xbG $ $ \xch $ $ \xbm (\xbS \xex X')
\xcc \xbG \xex X' $

$(\xbm *3)$ $ \xbm (\xbP' \xCK \xbS'') \xex X'' \xcc \xbm (\xbS)
\xex X''.$

A generalization to more than two factors is obvious.

One can also consider weakenings, e.g.,

$(S*1')$ $ \xbG' \xCK \xbS'' \xcc \xbS' \xCK \xbS'' $ is big iff $
\xbG' \xcc \xbS' $ is big.

\ed

\bfa

$\hspace{0.01em}$


\label{Fact Sin-Product-Small}

Let $(S*1)$ hold. Then:

$ \xbG' \xCK \xbG'' \xcc \xbS' \xCK \xbS'' $ is small iff $ \xbG'
\xcc \xbS' $ or $ \xbG'' \xcc \xbS'' $ is small.

\efa

\subparagraph{
Proof
}

$\hspace{0.01em}$


``$ \xci $'':

Suppose $ \xbG' \xcc \xbS' $ is small. Then $ \xbS' - \xbG' \xcc \xbS
' $ is big and
$(\xbS' - \xbG') \xCK \xbS'' \xcc \xbS' \xCK \xbS'' $ is big by
$(S*1).$ But
$(\xbG' \xCK \xbG'') \xcs ((\xbS' - \xbG') \xCK \xbS'')= \xCQ,$
so $ \xbG' \xCK \xbG'' \xcc \xbS' \xCK \xbS'' $ is small.

``$ \xch $'':

For the converse, suppose that neither $ \xbG' \xcc \xbS' $ nor $ \xbG
'' \xcc \xbS'' $ are small. Let
$A \xcc \xbS' \xCK \xbS'' $ be big, we show that $A \xcs (\xbG' \xCK
\xbG'') \xEd \xCQ.$
By $(S*1)$ there are $B' \xcc \xbS' $ and $B'' \xcc \xbS'' $ big, and
$B' \xCK B'' \xcc A.$ Then $B' \xcs \xbG' \xEd \xCQ,$
$B'' \xcs \xbG'' \xEd \xCQ,$ so there is $ \xBc x',x''  \xBe  \xbe (B' \xCK B''
)
\xcs (\xbG' \xCK \xbG'') \xcc A \xcs (\xbG' \xCK \xbG'').$

$ \xcz $
\\[3ex]

\bfa

$\hspace{0.01em}$


\label{Fact Sin-S-Mu}

If the filters over A are principal filters, generated by $ \xbm (A),$
i.e. $B \xcc A$ is big iff $ \xbm (A) \xcc B \xcc A$ for some $ \xbm (A)
\xcc A,$ then:

$(S*i)$ is equivalent to $(\xbm *i),$ $i=1,2,3.$

\efa

\subparagraph{
Proof
}

$\hspace{0.01em}$


(1)

``$ \xch $''

``$ \xcc $'': $ \xbm (\xbS') \xcc \xbS' $ and $ \xbm (\xbS'') \xcc
\xbS'' $ are big, so by $(S*1)$
$ \xbm (\xbS') \xCK \xbm (\xbS'') \xcc \xbS' \xCK \xbS'' $ is big,
so $ \xbm (\xbS' \xCK \xbS'') \xcc \xbm (\xbS') \xCK \xbm (\xbS''
).$

``$ \xcd $'': $ \xbm (\xbS' \xCK \xbS'') \xcc \xbS' \xCK \xbS'' $
is big $ \xch $ by $(S*1)$ there is $ \xbG' \xCK \xbG'' \xcc \xbm (\xbS
' \xCK \xbS'')$
and $ \xbG' \xcc \xbS',$ $ \xbG'' \xcc \xbS'' $ big $ \xch $ $ \xbm (
\xbS') \xcc \xbG',$ $ \xbm (\xbS'') \xcc \xbG'' $ $ \xch $
$ \xbm (\xbS') \xCK \xbm (\xbS'') \xcc \xbm (\xbS' \xCK \xbS''
).$

``$ \xci $''

Let $ \xbG' \xcc \xbS' $ be big, $ \xbG'' \xcc \xbS'' $ be big, $ \xbG
' \xCK \xbG'' \xcc \xbD,$ then $ \xbm (\xbS') \xcc \xbG',$
$ \xbm (\xbS'') \xcc \xbG'',$ so
by $(\xbm *1)$ $ \xbm (\xbS)= \xbm (\xbS') \xCK \xbm (\xbS'')
\xcc \xbG' \xCK \xbG'' \xcc \xbD,$ so $ \xbD $ is big.

Let $ \xbD \xcc \xbS $ be big, then by $(\xbm *1)$ $ \xbm (\xbS') \xCK
\xbm (\xbS'')= \xbm (\xbS) \xcc \xbD.$

(2)

``$ \xch $''

$ \xbm (\xbS) \xcc \xbG $ $ \xch $ $ \xbG \xcc \xbS $ big $ \xch $ by
$(S*2)$ $ \xbG \xex X' \xcc \xbS \xex X' $ big $ \xch $ $ \xbm (\xbS \xex
X') \xcc \xbG \xex X'.$

``$ \xci $''

$ \xbG \xcc \xbS $ big $ \xch $ $ \xbm (\xbS) \xcc \xbG $ $ \xch $ by $(
\xbm *2)$ $ \xbm (\xbS \xex X') \xcc \xbG \xex X' $ $ \xch $ $ \xbG \xex
X' \xcc \xbS \xex X' $ big.

(3)

``$ \xch $''

$ \xbm (\xbS) \xcc \xbS $ big $ \xch $ $ \xcE B \xcc \xbP' \xCK \xbS''
$ big s.t. $B \xex X'' \xcc \xbm (\xbS) \xex X'' $ by $(S*3),$ thus in
particular $ \xbm (\xbP' \xCK \xbS'') \xex X'' \xcc \xbm (\xbS) \xex
X''.$

``$ \xci $''

$A \xcc \xbS $ big $ \xch $ $ \xbm (\xbS) \xcc A.$ $ \xbm (\xbP' \xCK
\xbS'') \xcc \xbP' \xCK \xbS'' $ is big, and by $(\xbm *3)$
$ \xbm (\xbP' \xCK \xbS'') \xex X'' \xcc \xbm (\xbS) \xex X'' \xcc A
\xex X''.$

$ \xcz $
\\[3ex]

\paragraph{
Discussion
}

$\hspace{0.01em}$


\label{Section Discussion}

We compare these rules to probability defined size.

Let ``big'' be defined by
``more than $50 \xET $''. If $ \xbP X' $ and $ \xbP X'' $ have 3 elements
each, then subsets of $ \xbP X' $
or $ \xbP X'' $ of $card \xcg 2$
are big. But taking the product may give $4/9<1/2.$ So the product rule
``$big*big=big$'' will not hold there. One direction will hold, of course.

Next, we discuss the prerequisite $ \xbS = \xbS' \xCK \xbS''.$ Consider
the following
example:

\be

$\hspace{0.01em}$


\label{Example Sin-Prod-Size}

Take a language of 5 propositional variables, with $X':=\{a,b,c\},$ $X''
:=\{d,e\}.$
Consider the model set $ \xbS:=\{ \xCL a \xCL b \xCL cde,$ $-a-b-c-d \xCL
e\},$ i.e. of 8 models of
$ \xCf de$ and 2 models of $- \xCf d.$ The models of $ \xCf de$ are 8/10
of all elements of $ \xbS,$ so
it is reasonable to call them a big subset of $ \xbS.$ But its projection
on
$X'' $ is only 1/3 of $ \xbS''.$

So we have a potential $ \xCf decrease$ when going to the coordinates.

This shows that weakening the prerequisite about $ \xbS $ as done in
$(S*2)$ is not
innocent.

\ee

\br

$\hspace{0.01em}$


\label{Remark Sin-Boolean}

When we set small sets to 0, big sets to 1, we have the following boolean
rules
for filters:

(1) $0+0=0$

(2) $1+x=1$

(3) $-0=1,$ $-1=0$

(4) $0*x=0$

(5) $1*1=1$

There are no such rules for medium size sets, as the union of two medium
size sets may be big, but also stay medium.

Such multiplication rules capture the behaviour of Reiter defaults and of
defeasible inheritance.
\subsubsection{
Interpolation of the form $\xbf \xcl \xba \xcn \xbq $
}

\er

\bco

$\hspace{0.01em}$


\label{Corollary Sin-Interpolation-1}

We assume definability as shown in
Fact \ref{Fact Sin-Simplification-Definable} (page \pageref{Fact
Sin-Simplification-Definable}).

Interpolation of the form $ \xbf \xcl \xba \xcn \xbq $ exists, if

(1) both $(S*1)$ and $(S*2),$

or

(2) $(S*3)$ hold,

when $ \xbb \xcn \xbg $ is defined by:

$ \xbb \xcn \xbg: \xcj \xbm (\xbb)= \xbm (M(\xbb)) \xcc M(\xbg),$
and

$ \xbm (X)$ is the generator of the principal filter over $X.$

\eco

\subparagraph{
Proof
}

$\hspace{0.01em}$


Let $ \xbS:=M(\xbf),$ $ \xbG:=M(\xbq),$ $X' $ the set of variables
only in $ \xbf,$ so $ \xbG = \xbP' \xCK \xbG \xex X'',$
where $ \xbP' = \xbP X'.$
Set $ \xba:=Th(\xbP' \xCK \xbS''),$ where $ \xbS'' = \xbS \xex X''
.$
Note that variables only in $ \xbq $ are automatically taken care of,
as $ \xbS'' $ can be written as a product without mentioning them.

By prerequisite, $ \xbm (\xbS) \xcc \xbG,$ we have to show $ \xbm (
\xbP' \xCK \xbS'') \xcc \xbG.$

(1)

$ \xbm (\xbP' \xCK \xbS'')= \xbm (\xbP') \xCK \xbm (\xbS'')$ by
$(S*1)$ and
Fact \ref{Fact Sin-S-Mu} (page \pageref{Fact Sin-S-Mu})  (1).
By $ \xbm (\xbS) \xcc \xbG,$ $(S*2),$ and
Fact \ref{Fact Sin-S-Mu} (page \pageref{Fact Sin-S-Mu})  (2), $ \xbm (\xbS''
)= \xbm (\xbS \xex
X'') \xcc \xbG \xex X'',$ so
$ \xbm (\xbP' \xCK \xbS'')= \xbm (\xbP') \xCK \xbm (\xbS'') \xcc
\xbm (\xbP') \xCK \xbG \xex X'' \xcc \xbG.$

(2)

$ \xbm (\xbP' \xCK \xbS'') \xex X'' \xcc \xbm (\xbS) \xex X'' \xcc
\xbG \xex X'' $ by $(S*3)$ and
Fact \ref{Fact Sin-S-Mu} (page \pageref{Fact Sin-S-Mu})  (3).
So $ \xbm (\xbP' \xCK \xbS'') \xcc \xbP' \xCK (\xbm (\xbP' \xCK
\xbS'') \xex X'') \xcc \xbP' \xCK (\xbG \xex X'')= \xbG.$

$ \xcz $
\\[3ex]
\subsubsection{
Hamming relations
}

\bd

$\hspace{0.01em}$


\label{Definition Sin-Set-HR}

Call a relation $ \xec \xcc (\xbP X \xCK \xbP X) \xcv (\xbP' \xCK \xbP
') \xcv (\xbP'' \xCK \xbP'')$ (where $ \xbP' = \xbP X',$
$ \xbP'' = \xbP X'')$ a (generalized) Hamming relation iff

$ \xbs \xec \xbt $ $ \xcj $ $ \xbs' \xec \xbt' $ and $ \xbs'' \xec \xbt
'',$

where $ \xbs = \xbs' \xDO \xbs'',$ $ \xbt = \xbt' \xDO \xbt'' $ $(
\xDO $ concatenation).

Define $x \xeb y: \xcj x \xec y$ and $x \xEd y.$

Thus, $ \xbs \xeb \xbt $ iff $ \xbs \xec \xbt $ and $(\xbs' \xeb \xbt'
$ or $ \xbs'' \xeb \xbt'').$

Again, generalization to more than two components is straightforward.

\ed

\be

$\hspace{0.01em}$


\label{Example Sin-Int-Circum}

The circumscription relation is a Hamming relation with $ \xCN p \xck p$
and
$ \xcU \xCL q_{i} \xck \xcU \xCL q'_{i}$ iff $ \xcA i(\xCL q_{i} \xck
\xCL q'_{i}).$

\ee

\br

$\hspace{0.01em}$


\label{Remark Sin-Independence}

(1) The independence makes sense because the concept of models, and thus
the
usual
interpolation for classical logic relie on the independence of the
assignments.

(2) This corresponds to social choice for many independent dimensions.

(3) We can also consider such factorisation as an approximation:
we can do part of the reasoning independently.
\subsubsection{
Hamming relations and product size
}

\er

\bd

$\hspace{0.01em}$


\label{Definition Sin-Relation-to-Filter}

Given a relation $ \xec,$ define as usual a principal filter $ \xdf (X)$
generated by the
$ \xec -$minimal elements:

$ \xbm (X):=\{x \xbe X: \xCN \xcE x' \xeb x.x' \xbe X\},$

$ \xdf (X):=\{A \xcc X: \xbm (X) \xcc A\}.$

\ed

\bfa

$\hspace{0.01em}$


\label{Fact Sin-Component}

Let $ \xec $ be a smooth Hamming relation. Then
$(\xbm *2)$ holds, and thus $(S*2)$ by Fact \ref{Fact Sin-S-Mu} (page
\pageref{Fact Sin-S-Mu}),
(2).

\efa

\subparagraph{
Proof
}

$\hspace{0.01em}$


Suppose $ \xbm (\xbS) \xcc \xbG $ and $s \xbe \xbS \xex X' - \xbG \xex
X',$ we show $s \xce \xbm (\xbS \xex X').$

Let $ \xbs =s \xDO \xbs'' \xbe \xbS,$ then $ \xbs \xce \xbG,$ so $ \xbs
\xce \xbm (\xbS).$ So here is $ \xbr \xeb \xbs,$
$ \xbr \xbe \xbm (\xbS) \xcc \xbG $ by
smoothness. We have $ \xbr' \xec s.$ $ \xbr' =s$ cannot be, as $ \xbr'
\xbe \xbG \xex X',$ and
$s \xce \xbG \xex X'.$ So $ \xbr' \xeb s,$ and $s \xce \xbm (\xbS \xex
X').$

$ \xcz $
\\[3ex]

\be

$\hspace{0.01em}$


\label{Example Sin-Component-Inverse}

Even for smooth Hamming relations, the converse of $(\xbm *2)$ is not
necessarily
true:

Let $ \xbs' \xeb \xbt',$ $ \xbt'' \xeb \xbs'',$ $ \xbs:=\{ \xbs,
\xbt \},$ then $ \xbm (\xbS)= \xbS,$ but $ \xbm (\xbS')=\{ \xbs'
\},$
$ \xbm (\xbS'')=\{ \xbt'' \},$ so $ \xbm (\xbS) \xEd \xbm (\xbS')
\xCK \xbm (\xbS'').$

We need the additional assumption that $ \xbm (\xbS') \xCK \xbm (\xbS
'') \xcc \xbS,$ see
Fact \ref{Fact Sin-Product} (page \pageref{Fact Sin-Product})  (1).

\ee

\bfa

$\hspace{0.01em}$


\label{Fact Sin-Product}

Let again $ \xbS':= \xbS \xex X',$ $ \xbS'':= \xbS \xex X''.$

(1) Let $ \xec $ be a smooth Hamming relation. Then:

$ \xbm (\xbS') \xCK \xbm (\xbS'') \xcc \xbS $ $ \xch $ $ \xbm (\xbS
)= \xbm (\xbS') \xCK \xbm (\xbS'').$

(Here not necessarily $ \xbS = \xbS' \xCK \xbS''.)$

(2) Let $ \xec $ be a Hamming relation, and $ \xbS = \xbS' \xCK \xbS''
.$ Then
$(\xbm *1)$ holds, and thus, by Fact \ref{Fact Sin-S-Mu} (page \pageref{Fact
Sin-S-Mu}), (1)
$(S*1).$

\efa

\subparagraph{
Proof
}

$\hspace{0.01em}$


(1)

``$ \xcd $'': Let $ \xbs' \xbe \xbm (\xbS'),$ $ \xbs'' \xbe \xbm (
\xbS'').$ By prerequisite, $ \xbs' \xDO \xbs'' \xbe \xbS.$ Suppose
$ \xbt \xeb \xbs' \xDO \xbs'',$ then $ \xbt' \xeb \xbs' $ or $ \xbt
'' \xeb \xbs'',$ contradiction.

``$ \xcc $'': Let $ \xbs \xbe \xbm (\xbS),$ suppose $ \xbs' \xce \xbm
(\xbS')$ or $ \xbs'' \xce \xbm (\xbS'').$
So there are $ \xbt' \xec \xbs',$ $ \xbt'' \xec \xbs'' $ with $ \xbt
' \xbe \xbm (\xbS'),$ $ \xbt'' \xbe \xbm (\xbS'')$ by
smoothness. Moreoever, $ \xbt' \xeb \xbs' $ or $ \xbt'' \xeb \xbs''.$
By prerequisite $ \xbt' \xDO \xbt'' \xbe \xbS,$
and $ \xbt' \xDO \xbt'' \xeb \xbs,$ so $ \xbs \xce \xbm (\xbS).$

(2)

``$ \xcd $'': As in (1), the prerequisite holds trivially.

``$ \xcc $'': As in (1), but we do not need $ \xbt' \xbe \xbm (\xbS'
),$ $ \xbt'' \xbe \xbm (\xbS''),$ as $ \xbt' \xDO \xbt'' $ will
be in $ \xbS $ trivially. So smoothness is not needed.

$ \xcz $
\\[3ex]

We now show the main property for $ \xcl \xDO \xcn $ interpolation in a
direct proof.
Note that interpolation treats the two components of the product
differently, which also shows that our requirements are sufficient, but
not necessary - see
Section \ref{Section Sin-Context} (page \pageref{Section Sin-Context}).

\bfa

$\hspace{0.01em}$


\label{Fact Sin-Int-Prod}

Let $ \xec $ be a smooth Hamming relation, then
$(\xbm *3)$ holds, and thus by Fact \ref{Fact Sin-S-Mu} (page \pageref{Fact
Sin-S-Mu})  (3) $(S*3).$

\efa

\subparagraph{
Proof
}

$\hspace{0.01em}$


Let again $ \xbP' = \xbP X',$ $ \xbS'' = \xbS \xex X''.$

Let $ \xbD:= \xbP' \xCK \xbS'',$ $ \xbs'' \xbe \xbm (\xbD) \xex X''
,$ $ \xbs = \xbs' \xDO \xbs'' \xbe \xbm (\xbD)$ for some $ \xbs' \xbe
\xbP'.$
Suppose $ \xbs'' \xce \xbm (\xbS) \xex X''.$ There cannot be any $
\xbt \xeb \xbs,$ $ \xbt \xbe \xbS,$ by $ \xbS \xcc \xbD.$
So $ \xbs \xce \xbS,$ but $ \xbs'' \xbe \xbS'',$ so there is $ \xbt
\xbe \xbS $ $ \xbt'' = \xbs''.$ As $ \xbt $ is not minimal,
there must be minimal $ \xbr = \xbr' \xDO \xbr'' \xeb \xbt,$ $ \xbr
\xbe \xbS $ by smoothness. As $ \xbr $ is minimal,
$ \xbr'' \xEd \xbs'',$ and as
$ \xbr \xeb \xbt,$ $ \xbr'' \xeb \xbs''.$ By prerequisite $ \xbs'
\xDO \xbr'' \xbe \xbD,$ and $ \xbs' \xDO \xbr'' \xeb \xbs,$
contradiction.

Note that smoothness is essential. Otherwise, there might be an infinite
descending chain $ \xbt_{i}$ below $ \xbt,$ all with $ \xbt''_{i}= \xbs
'',$ but none below $ \xbs.$

$ \xcz $
\\[3ex]

\bco

$\hspace{0.01em}$


\label{Corollary Sin-Interpolation-2}

Interpolation in the form $ \xbf \xcl \xba \xcn \xbq $ exists, when $ \xcn
$ is defined by a smooth
Hamming relation.

\eco

\subparagraph{
Proof
}

$\hspace{0.01em}$


We give two proofs:

(1)

By Fact \ref{Fact Sin-Component} (page \pageref{Fact Sin-Component})  and
Fact \ref{Fact Sin-Product} (page \pageref{Fact Sin-Product})  (2) $(S*1)$ and
$(S*2)$ hold.
Thus, by Corollary \ref{Corollary Sin-Interpolation-1} (page \pageref{Corollary
Sin-Interpolation-1})  (1), interpolation
exists.

(2)

By Fact \ref{Fact Sin-Int-Prod} (page \pageref{Fact Sin-Int-Prod}), $(S*3)$
holds, so by
Corollary \ref{Corollary Sin-Interpolation-1} (page \pageref{Corollary
Sin-Interpolation-1})  (2), interpolation exists.

$ \xcz $
\\[3ex]

\bfa

$\hspace{0.01em}$


\label{Fact Sin-Representation}

$(\xbm *1)$ and $(\xbm *2)$ and the usual axioms for smooth relations
characterize smooth Hamming relations.

\efa

\subparagraph{
Proof
}

$\hspace{0.01em}$


Define as usual $ \xbs \xeb \xbt: \xcj \xbt \xce \xbm (\{ \xbs, \xbt
\}).$

We have to show:

$ \xbs \xeb \xbt $ iff $ \xbs' \xec \xbt' $ and $ \xbs'' \xec \xbt'' $
and $(\xbs' \xeb \xbt' $ or $ \xbs'' \xeb \xbt'').$

``$ \xci $'':

Suppose $ \xbs' \xeb \xbt' $ and $ \xbs'' \xec \xbt''.$ Then $ \xbm
(\{ \xbs', \xbt' \})=\{ \xbs' \},$ and
$ \xbm (\{ \xbs'', \xbt'' \})=\{ \xbs'' \}$ (either $ \xbs'' \xeb
\xbt'' $ or $ \xbs'' = \xbt'',$ so in both cases
$ \xbm (\{ \xbs'', \xbt'' \})=\{ \xbs'' \}).$ As $ \xbt' \xce \xbm
(\{ \xbs', \xbt' \}),$
$ \xbt \xce \xbm (\{ \xbs', \xbt' \} \xCK \{ \xbs'', \xbt'' \})=
\xbm (\{ \xbs', \xbt' \}) \xCK \xbm (\{ \xbs'', \xbt'' \})$ (by $(
\xbm *1))$
$=\{ \xbs' \} \xCK \{ \xbs'' \}=\{ \xbs \},$
by $(\xbm *1),$ so by smoothness $ \xbs \xeb \xbt.$

``$ \xch $'':

Conversely, if $ \xbs \xeb \xbt,$ so $ \xbG:=\{ \xbs \}= \xbm (\xbS)$
for $ \xbS:=\{ \xbs, \xbt \},$ so by $(\xbm *2)$
$ \xbm (\xbS \xex X')= \xbm (\{ \xbs', \xbt' \}) \xcc \xbG \xex X'
=\{ \xbs' \},$ so $ \xbs' \xec \xbt',$
analogously $ \xbm (\xbS \xex X'')= \xbm (\{ \xbs'', \xbt'' \}) \xcc
\xbG \xex X'' =\{ \xbs'' \},$ so $ \xbs'' \xec \xbt'',$ but both
cannot be equal.

$ \xcz $
\\[3ex]
\subsection{
Context and structure
}

\label{Section Sin-Context}

We take the importance of condition $(\xbm *3)$ (or $(S*3))$
as occasion for a broader remark.

 \xEh

 \xDH
This condition points to a weakening of the Hamming condition:

Adding new ``branches'' in $X' $ will not give new minimal elements in $X''
,$
but may destroy other minimal elements in $X''.$ This can be achieved
by a sort of semi-rankedness:
If $ \xbr $ and $ \xbs $ are different only in the $X' -$part, then $ \xbt
\xeb \xbr $ iff $ \xbt \xeb \xbs,$
but not necessarily $ \xbr \xeb \xbt $ iff $ \xbs \xeb \xbt.$

 \xDH
In more abstract terms:

When we separate support from attack (support: a branch $ \xbs' $ in $X'
$ supports
a continuation $ \xbs'' $ in $X'' $ iff $ \xbs \xDO \xbs'' $ is minimal,
i.e. not attacked, attack:
a branch $ \xbt $ in $X' $ attacks a continuation $ \xbs'' $ in $X'' $
iff it prevents all $ \xbs \xDO \xbs'' $
to be minimal), we see that new branches will not support any new
continuations,
but may well attack continuations.

More radically, we can consider paths $ \xbs'' $ as positive information,
$ \xbs' $ as potentialy negative information. Thus, $ \xbP' $ gives
maximal
negative information, and thus smallest set of accepted models.

 \xDH
We can interpret this as follows:
$X'' $ determines the base set.
$X' $ is the context. This determines the choice (subset of the base set).
We compare to preferential structures:
In preferential structures, $ \xeb $ is not part of the
language either, it is context. And we have the same behaviour as
shown in the
fundamental property of preferential structures:
the bigger the set, the more attacks are possible.

 \xDH
The concept of size looks only at the result of support and attack, so it
is necessarily somewhat coarse. Future research should also investigate
both
concepts separately.

 \xEj

We broaden this.

Following a tradition begun by Kripke, one has added structure to the set
of
classical models, reachability, preference, etc. Perhaps one should
emphasize a
more abstract approach,
as done by one the authors e.g. in  \cite{Sch92}, and elaborated in
 \cite{Sch04}, see in particular the distinction between structural
and algebraic
semantics in the latter. Our suggestion is to separate structure from
logic in
the semantics, and to treat what we called context above by a separate
``machinery''. So, given a set $X$ of models, we have some abstract
function $f,$ which chooses the models where the consequences hold,
$f(X).$

Now, we can put into this ``machinery'' whatever we want.

The abstract properties of preferential or modal structures are well
known.

But we can also investigate non-static $f,$ where $f$ changes in function
of what we already did - ``reacting'' to the past.

We can look at usual properties of $f,$ complexity,
generation by some simple structure like a special
machine, etc.

So we advocate the separation of usual, classical semantics,
from the additional properties, which are treated
``outside''.
It might be interesting to forget altogether about logic,
classify those functions or more complicated devices which
correspond to some logical property, and investigate them and
their properties.
\section{
Hamming distances and revision
}

$ \xCO $
\index{Hamming distance}

\bd

$\hspace{0.01em}$


\label{Definition Hamming-Distance}

Given $x,y \xbe \xbS,$ a set of sequences over an index set $I,$
the Hamming distance \index{Hamming distance}
comes in two flavours:

$d_{s}(x,y)$ \index{$d_{s}(x,y)$}  $:=\{i \xbe I:x(i) \xEd y(i)\},$ the set
variant,

$d_{c}(x,y)$ \index{$d_{c}(x,y)$}  $:=card(d_{s}(x,y)),$ the counting variant.

We define $d_{s}(x,y) \xck d_{s}(x',y')$ iff $d_{s}(x,y) \xcc d_{s}(x'
,y'),$

thus, $s-$distances are not always comparabel.

There are straightforward generalizations of the counting variant:

We can also give different importance to different $i$ in the counting
variant, so e.g.,
$d_{c}(\xBc x,x'  \xBe, \xBc y,y'  \xBe)$ might be 1 if $x \xEd y$ and $x'
=y',$ but 2 if
$x=y$ and $x' \xEd y'.$

If the $x \xbe \xbS $ may have more than 2 different values, then a
varying
individual distance may also reflect to the distances in $ \xbS.$
So, if $d(x(i),x' (i))<d(x(i),x'' (i)),$ then (the rest being equal),
we may have $d(x,x')<d(x,x'').$

\ed

$ \xCO $

The (for us) essential property of the set variant Hamming distance is
that we
cannot compensate, i.e. a difference in $x$ cannot be compensated by
equality
in $x' $ and $x'',$ as in the counting variant.

\bd

$\hspace{0.01em}$


\label{Definition Sin-Set-HD}

Call a function $d:(\xbP X \xCK \xbP X) \xcv (\xbP' \xCK \xbP') \xcv
(\xbP'' \xCK \xbP'') \xcp Z$

(where $ \xbP' = \xbP X',$ $ \xbP'' = \xbP X'')$

a (generalized) Hamming distance iff

(1) there is a total order $ \xck $ on $Z$

(2) $d(\xbs_{1}, \xbt_{1}) \xck d(\xbs_{2}, \xbt_{2})$ iff
$d(\xbs'_{1}, \xbt'_{1}) \xck d(\xbs'_{2}, \xbt'_{2})$ and $d(\xbs
''_{1}, \xbt''_{1}) \xck d(\xbs''_{2}, \xbt''_{2}).$

(In the strict part, at least one has to be strict on the right hand
side.)

\ed

\bd

$\hspace{0.01em}$


\label{Definition Sin-Bar}

Given a distance $d,$ define for two sets $X,Y$

$X \xfA Y:=\{y \xbe Y: \xcE x \xbe X(\xCN \xcE x' \xbe X,y' \xbe Y.d(x'
,y')<d(x,y))\}.$

We assume that $X \xfA Y \xEd \xCQ $ if $X,Y \xEd \xCQ.$

\ed

We have results analogous to the relation case.

\bfa

$\hspace{0.01em}$


\label{Fact Sin-HD-Product}

Let $ \xfA $ be defined by a Hamming distance, then:

(1) $(\xbS_{1}' \xCK \xbS_{1}'') \xfA (\xbS_{2}' \xCK \xbS_{2}'')=(
\xbS_{1}' \xfA \xbS_{2}') \xCK (\xbS_{1}'' \xfA \xbS_{2}'')$

(2) $(\xbS_{1}' \xfA \xbS_{2}') \xCK (\xbS_{1}'' \xfA \xbS_{2}'') \xcc
\xbS_{2}$ and $(\xbS_{2}' \xfA \xbS_{1}') \xCK (\xbS_{2}'' \xfA
\xbS_{1}'') \xcc \xbS_{1}$ $ \xch $
$(\xbS_{1}) \xfA (\xbS_{2})=(\xbS_{1}' \xfA \xbS_{2}') \xCK (
\xbS_{1}'' \xfA \xbS_{2}''),$ if the distance is symmetric

(where $ \xbS_{i}$ here is not necessarily $ \xbS_{i}' \xCK \xbS_{i}'',$
etc.).

\efa

\subparagraph{
Proof
}

$\hspace{0.01em}$


(1) and (2)

``$ \xcc $'':

Suppose $d(\xbs, \xbt)$ is minimal for $ \xbs \xbe \xbS_{1},$ $ \xbt
\xbe \xbS_{2}.$ If there is
$ \xba' \xbe \xbS_{1}',$ $ \xbb' \xbe \xbS_{2}' $ s.t. $d(\xba',
\xbb')<d(\xbs', \xbt'),$ then
$d(\xba' \xDO \xbs'', \xbb' \xDO \xbt'')<d(\xbs, \xbt),$ so $d(
\xbs', \xbt')$ and $d(\xbs'', \xbt'')$ have to be
minimal.

For the argument to go through, we need $ \xba' \xDO \xbs'',$ $ \xbb'
\xDO \xbt'' $ to be
in the sets considered. If $ \xbS_{i}= \xbS_{i}' \xCK \xbS_{i}'',$ this
is trivially satisfied.
Otherwise, we use the usual minimality assumption of revision, this
replaces smoothness, and add the supplementary conditions.

``$ \xcd $'':

For the converse, suppose $d(\xbs', \xbt')$ and $d(\xbs'', \xbt''
)$ are minimal, but
$d(\xbs, \xbt)$ is not, so $d(\xba, \xbb)<d(\xbs, \xbt)$ for some
$ \xba, \xbb,$ then
$d(\xba', \xbb')<d(\xbs', \xbt')$ or $d(\xba'', \xbb'')<d(
\xbs'', \xbt''),$ contradiction.

$ \xcz $
\\[3ex]

\bco

$\hspace{0.01em}$


\label{Corollary Sin-Parikh-Revision}

By Fact \ref{Fact Sin-HD-Product} (page \pageref{Fact Sin-HD-Product}), Hamming
distances
generate decomposable revision operators a la Parikh,
see  \cite{Par96},
also in the generalized form of variable $K$ and $ \xbf.$

\eco

\bfa

$\hspace{0.01em}$


\label{Fact Sin-HD-Part}

Let $ \xfA $ be defined by a Hamming distance, then:

$ \xbP \xfA \xbS \xcc \xbG $ $ \xch $ $ \xbP' \xfA (\xbS \xex X') \xcc
\xbG \xex X'.$

\efa

\subparagraph{
Proof
}

$\hspace{0.01em}$


Let $t \xbe \xbS \xex X' - \xbG \xex X',$ we show $t \xce \xbP' \xfA (
\xbS \xex X').$ Let $ \xbt \xbe \xbS $ be s.t.
$ \xbt' =t,$ then $ \xbt \xce \xbG $ (otherwise $t \xbe \xbG \xex X'),$
so $ \xbt \xce \xbP \xfA \xbS,$ so there is $ \xba \xbe \xbP,$
$ \xbb \xbe \xbS,$ $d(\xba, \xbb)<d(\xbs, \xbt)$ for all $ \xbs
\xbe \xbP.$ If $d(\xbs', \xbt')$ were minimal
for some $ \xbs,$ then we would consider $ \xbs' \xDO \xba'',$ $ \xbt
' \xDO \xbb'',$ and have
$ \xbt' \xDO \xbb'' \xbe \xbP \xfA \xbS,$ so $ \xbt' \xDO \xbb'' \xbe
\xbG,$ and $t \xbe \xbG \xex X',$ contradiction.

$ \xcz $
\\[3ex]
\subsection{
Discussion of representation
}

It would be nice to have a representation result like the one for
Hamming relations. But this is impossible, for the following reason:

In constructing the representing distance from revision results, we made
arbitrary choices (see the proofs in  \cite{LMS01} or
 \cite{Sch04}). I.e., we choose sometimes
arbitrarily $d(x,y) \xck d(x',y'),$ when we do not have enough
information
to decide. (This is an example of the fact that the problem of
``losing ignorance'' should not be underestimated, see e.g.
 \cite{GS08f}.)
As we do not follow the same procedure for all cases, there
is no guarantee that the different representations will fit together.

Of course, it might be possible to come to a uniform choice, and one
could then attempt a representation result. This is left as an open
problem.

\end{document}